\documentclass[a4paper,reqno,10pt,oneside]{amsart}

\usepackage{amsmath,amsthm,amsfonts,amssymb}
\usepackage{indentfirst}
\usepackage{url}
\usepackage{hyperref}
\usepackage{graphicx}
\usepackage{color}
\usepackage{a4wide}
\usepackage{setspace}
\usepackage{comment}

\usepackage{enumerate}

\theoremstyle{plain}
\newtheorem{theorem}{Theorem}[section]

\theoremstyle{definition}
\newtheorem{defn}[theorem]{Definition}
\newtheorem{exa}[theorem]{Example}
\newtheorem{obs}[theorem]{Remark}

\newcommand{\xx}{\left(\frac{x}{x_{2,0}}\right)}

\begin{document}

	\title[A modified Maki-Thompson model with directed inter-group interactions]
	{Asymptotic behavior for a modified Maki-Thompson model\\ with directed inter-group interactions}
	\author[Carolina Grejo]{Carolina Grejo}
	\author[Pablo M. Rodriguez]{Pablo M. Rodriguez}
	\address[C. Grejo and P. M. Rodriguez]{Instituto de Ci\^encias Matem\'aticas e de Computa\c{c}\~ao, Universidade de S\~ao Paulo. Caixa Postal 668, 13560-970 S\~ao Carlos, SP, Brazil.
}
	\email{pablor@icmc.usp.br}
	\email{carolina@ime.usp.br}

	\thanks{}

	\date{\today}
	
	\begin{abstract} 
	In this work we propose a new extension for the Maki-Thompson rumor model which incorporates inter-group directed contacts. The model is defined on an homogeneously mixing population where the existence of two differentiated groups of individuals is assumed. While individuals of one group have an active role in the spreading process, individuals of the other group only contribute in stifling the rumor provided they would contacted. For this model we measure the impact of dissemination by studying the remaining proportion of ignorants of both groups at the end of the process. In addition we discuss some examples and possible applications.      
	\end{abstract}
	
	\maketitle
	
	
	\section{Introduction}
	
The mathematical modeling of rumor spreading has evidenced an increasing interest during the last years. An important part of the recent research has been motivated by the search for generalizations of two basic and well-accepted theoretical models; namely the Daley-Kendall and the Maki-Thompson models. The former was proposed in the midd 60's as an alternative to the well-known susceptible-infected-removed epidemic model (see \cite{daley_nature,kendall}), and the latter appeared in the literature as its simplification (see \cite{maki-1973}). The basic versions of these mathematical models intend to illustrate in a simple way how a rumor can spread through a homogeneously mixing population. With this application in mind, the Maki-Thomspon model is formulated by assuming that the population is subdivided into three classes of individuals; namely, ignorants, spreaders, and stiflers. Spreaders try to tell the rumor to any other individual at the instants of independent Poisson processes. If the contacted individual is an ignorant, then he/she becomes a spreader; otherwise the spreader who is trying to tell the rumor becomes a stifler. So stiflers represent individuals who know about the rumor but they are no longer being part of the spreading process. The quantities of people belonging to each class at any time can be described by a continuous-time Markov chain. On a finite population, what happens is that the process eventually ends and, when that occurs, the interest is to know the remaining proportion of people who never hear about the rumor. This is one way to measure the range of the spreading process. Indeed, the first rigorous results for these models are limit theorems for the remaining proportion of ignorants when the process ends, and it has been proved that such a proportion approximates to $20\%$ of the population, see \cite{sudbury,watson}. We refer the reader to \cite{dg,Lebensztayn-JMAA2015,lebensztayn/machado/rodriguez/2011b,lebensztayn/machado/rodriguez/2011a} for a review of results and generalizations of these models. Also, many modifications of these models have been considered assuming that the population is not necessarily homogeneous nor totally mixing, see \cite{EAP,arruda,raey,MNP-PRE2004,moreno-PhysA2007} and the references therein. The results obtained in that direction rely mainly on mean-field approximations, computational simulations, or partial qualitative results. 

The assumption of the population be homogeneously mixing allows to obtain many rigorous results, even if we incorporate more complex interactions between the classes of individuals. The reader can obtain a good illustration of this in the general models proposed in \cite{lebensztayn/machado/rodriguez/2011b,lebensztayn/machado/rodriguez/2011a}. In contrast, we suggest the reader to see \cite{arruda} for an overview of how to deal with complex interactions on an heterogeneous population. 

The motivation of this work is to propose and analyze a general Maki-Thompson model which incorporates directed inter-group interactions. As far as we know no rigorous results exist for this type of model and therefore it may be seen as a contribution to increase the familiy of general models like those considered recently by \cite{lebensztayn/machado/rodriguez/2011b,lebensztayn/machado/rodriguez/2011a}. In our model we assume the existence of two groups of individuals, say $A$-individuals and $B$-individuals. On the one hand, $A$-individuals are subdivided as usual into ignorants, spreaders and stiflers, and spreaders can transmit the rumor to any other individual of the population (but not necessarily at the same rate to $A$-individuals than to $B$-individuals). Inside the group of $A$-individuals the rumor spreads following the rules of the Maki-Thompson model. On the other hand, $B$-individuals are subdivided only into ignorants and stiflers. They role is passive with respect to spreading but each $B$-individual contributes in stifling the rumor if he/she would contacted at least twice by an $A$-individual who is a spreader. This type of dynamic is inspired by the GBN-model proposed in \cite{Brooks} as a mathematical model of negative rumor spread in the context of conflicting groups. Our goal is to determine the asymptotic proportions of ignorants at the end of the process as the population size grows to infinity.

The paper is organized as follows. Section \ref{section:main} includes the formal definition of the model and the main results as well as some examples. Our results are a Weak Law of Large Numbers and a Central Limit Theorem. Section \ref{section:proofs} is devoted to a discussion of the main steps and arguments to prove the theorems. Our arguments rely on the application of convergence results for density dependent Markov chains. This has been a useful tool to deal with this type of models, see for example \cite{lebensztayn/machado/rodriguez/2011b,lebensztayn/machado/rodriguez/2011a}, and it is an alternative to the pgf method and the Laplace transform
presented in \cite{dg,gani-Env2000,pearce}.

	\section{The model and main results}\label{section:main}
	We consider a population with $N$ individuals, of which $N_1$ belong to a first group (called $A$-individuals) and $N_2$ belongs to a second group (called $B$-individuals). We assume $N_1:= \theta\, N$ and $N_2:= (1-\theta)\, N$ where $\theta \in [0,1]$. Let $X_{1}^N(t), X_{2}^N(t), Y_{1}^N(t), Z^N(t)$ be the number of $A$-ignorants, $B$-ignorants, $A$-spreaders, and stiflers at time $t$, for $t\geq 0$, and let us assume $X_{1}^N(t)+X_{2}^N(t)+Y_{1}^N(t)+Z^N(t)=N$ for any $t\geq 0$; i.e. we consider a finite closed population. In addition, we assume that the following limits exists:
\begin{equation}\label{eq:limits_0}
\lim_{N\to \infty} \frac{X_{1}^N(t)}{N} = x_{1,0}>0,\;\;\;\; \lim_{N\to \infty} \frac{X_{2}^N(t)}{N} = x_{2,0}>0,\;\;\;\;\lim_{N\to \infty} \frac{Y_{1}^N(t)}{N} = y_{1,0},\;\;\;\;\lim_{N\to \infty} \frac{Z^N(t)}{N} = z_{0}.
\end{equation}	
Note that $x_{1,0},x_{2,0},y_{1,0},z_{0}\in [0,1]$ and $x_{1,0} + x_{2,0} + y_{1,0} + z_{0}=1$. We define the general Maki-Thompson model with parameters $\theta,\lambda,\alpha,p$ and initial conditions $x_{1,0},x_{2,0},y_{1,0},z_{0}$ as the continuous-time Markov chain $\left\{\left(X^N_1(t),X^N_2(t),Y^N_1(t)\right)\right\}_{t\in[0,\infty)}$, with transitions and rates given by:
\begin{equation}\label{eq:rates_general}
\begin{array}{cc}
\mbox{transition} \quad &\mbox{rate} \\[0.2cm]
(-1,0,1)       \quad & \lambda \,p\, X_{1}\,Y_1, \\[0.2cm]
(0,-1,0)     \quad & \alpha\, Y_1\, X_2, \\[0.2cm]
(0,0,-1)      \quad & \alpha \,Y_1\,(N_2-X_2) + \lambda\, Y_1\,(N_1-X_1),\\[0.2cm]
(-1,0,0)      \quad & \lambda\, (1-p)\, X_1\, Y_1, 
\end{array}
\end{equation}

\smallskip
\noindent
where $\lambda>0, \alpha >0$ and $p\in [0,1]$. In what follows we use the notation $\mathcal{X}_{GMT}(\theta,\lambda,\alpha,p;x_{1,0},x_{2,0},y_{1,0})$ for this Markov chain. In words, the model incorporates four possible transitions summarized in \eqref{eq:rates_general}. $\lambda$ is the rate at which an $A$-spreader contact any other $A$-individual, while $\alpha$ is the rate at which the contact is done with a $B$-individual. Once an $A$-ignorant is contacted by an $A$-spreader, the first becomes an $A$-spreader with probability $p$ or a stifler otherwise. Finally, any $A$-spreader interacting with another $A$-spreader, or a stifler, becomes a stifler. We point out that the basic Maki-Thompson model is obtained by considering $\theta = \lambda = p =1$, and $x_{1,0}=1, y_{1,0}=0$ (these initial conditions come from $X_1^N(0)=N-1$ and $Y_1^N(0)=1$ in the basic version).

\begin{obs}
Thus defined, our model increases the range of generalizations obtained by \cite{lebensztayn/machado/rodriguez/2011b,lebensztayn/machado/rodriguez/2011a}. The novelty lies in the incorporation of intergroup interactions. In order to illustrate the applicability of the model we propose two scenarios:
\begin{enumerate}[(i)]
\item {\it Negative rumor spreading in the context of conflicting groups.} Motivated by the GBN-model introduced in \cite{Brooks} we can assume that each individual is identified with a fixed group membership, either non-targeted ($A$-individuals) or maligned ($B$-individuals), with respect to the rumor target. As an example of this scenario let us cite the one proposed by \cite{Brooks}: in a discussion of a negative rumor about university faculty, professors would be considered members of the maligned group, while students would be part of the non-targeted group. In this case, a natural assumption is to consider $\lambda > \alpha$.

\smallskip
\item {\it Rumor propagation with the presence of identified observers.} This situation is related to experimental studies of rumor spreading. We can assume the existence of two separate groups of individuals involved in the transmission of rumors. This example is motivated in the experimental analysis performed by \cite{walker}, where comparable dread and wish rumors were planted concurrently in a college community.  On one hand, one group receives the rumor plants ($A$-individuals) and, on the other hand, a second group serves as rumor reporters ($B$-individuals). It can be assumed that those who receive the plants would pass the rumor on to others, with the planted rumors eventually reaching some individuals in the group of rumor reporters. If rumor reporters are oriented to not take part of the spreading procedure we obtain our proposed model. In this case the natural assumption is $\lambda=\alpha$.

\end{enumerate}
\end{obs}

\smallskip
As usual we define the absorption time of the process as ${\tau}^{(N)}:=\inf\{t\geq 0: Y_1^{N}(t)=0\}$, and we study the remaining proportions of ignorants,  $X_{1}^N({\tau}^{(N)})/N$ and $X_{2}^N({\tau}^{(N)})/N$, when the population size growths to infinity. The following definition will be useful to localize those limiting values.

\smallskip
\smallskip

\begin{defn}\label{def:xinf2}
 Let $f:(0,x_{2,0}]\to \mathbb{R}$ be the function given by
	\begin{equation}
	f(x)=y_{1,0}+ (1+p)x_{1,0}\left[1- \xx^{\lambda/\alpha}\right] +(x_{2,0}-x)+\left[\frac{(\lambda - \alpha)\theta+\alpha}{\alpha}\right]\ln\xx.
\end{equation}

\smallskip
\noindent
Define $x_{2,\infty}:=x_{2,\infty}(\lambda,\alpha,\theta,x_{1,0},x_{2,0},y_{1,0})$ as the unique root of $f$ in $(0,x_{2,0}]$ satisfying $f'(x)\geq 0$.  In addition, let  $x_{1,\infty}:=x_{1,0}\left(x_{2,{\infty}}/x_{2,0}\right) ^{\lambda/\alpha}$. 
\end{defn}

\begin{figure}[h!]

\includegraphics[scale=.3]{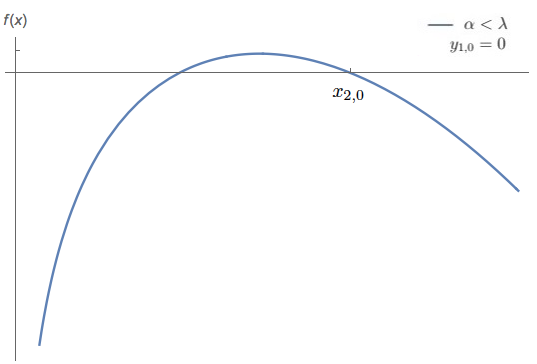}\hspace{1cm}
\includegraphics[scale=.3]{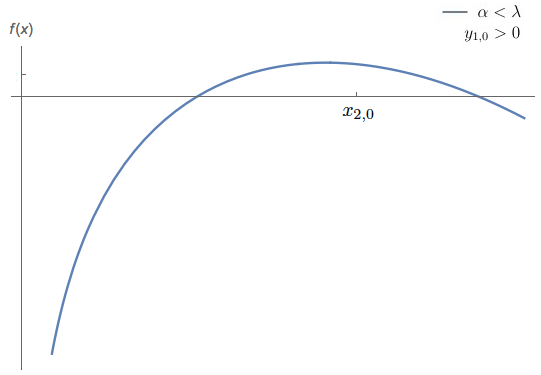}

\vspace{.8cm}
\includegraphics[scale=.3]{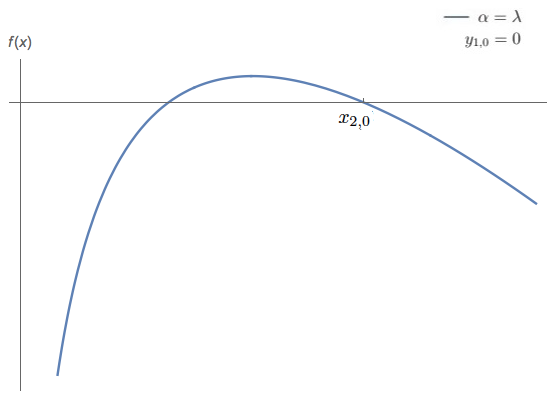}\hspace{1cm}
\includegraphics[scale=.3]{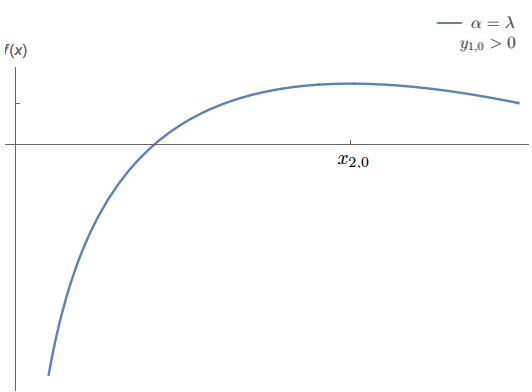}

\vspace{.8cm}
\includegraphics[scale=.3]{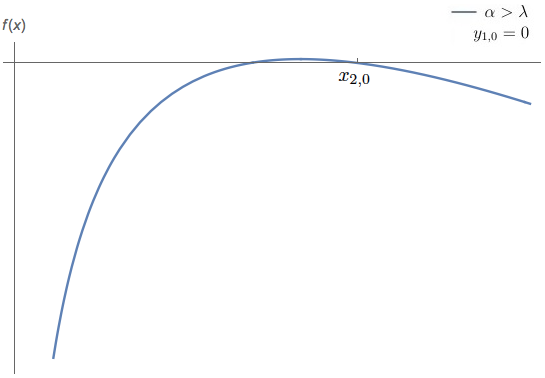}\hspace{1cm}
\includegraphics[scale=.3]{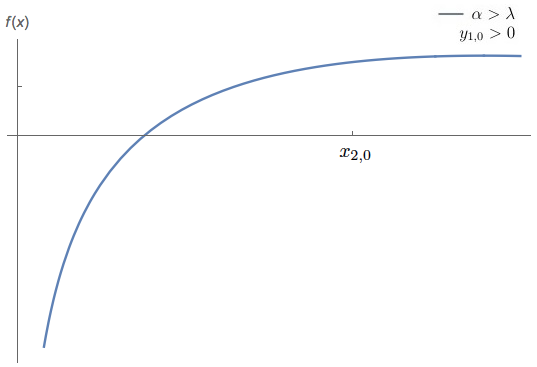}

\caption{Illustration of the behavior of the function $f$ according to $\alpha$, $\lambda$ and $y_{1,0}$.}\label{fig:function}

\end{figure}

\begin{obs}
As we shall see later the values $x_{1,\infty}$ and $x_{2,\infty}$ are the asymptotic remaining proportions of $A$-ignorants and $B$-ignorants respectively. It is not difficult to see that these quantities are well-defined. To see that, it is enough to verify that $f$ is continuous in $(0,x_{2,0}]$, $\lim_{x\searrow 0}f(x)=-\infty$ and $f(x_{2,0})=y_{1,0}\geq 0$. Moreover, $f$ is a concave function provided $\lambda \geq \alpha$. The case $\lambda<\alpha$ requires a bit more calculations. Indeed, let 

$$a(\lambda,\alpha,p,x_{1,0},x_{2,0}):=\frac{\alpha-\alpha x_{2,0}-\lambda x_{1,0}(1+p)}{\alpha - \lambda}$$

\smallskip
\noindent
and note that $f'(x_{2,0})<0$ if $\theta \in (a(\lambda,\alpha,p,x_{1,0},x_{2,0}), 1]$. On the other hand, let

$$b(\lambda,\alpha,p,x_{1,0}):=\frac{{\alpha}^2-\lambda x_{1,0}(1+p)(\alpha-\lambda)}{\alpha(\alpha - \lambda)},$$

\smallskip
\noindent
and 
$${\bar{x}}:=x_{2,0}\left[ \frac{(\alpha +(\lambda - \alpha) \theta)\alpha}{x_{1,0}\lambda(\alpha-\lambda)(1+p)}\right]^{\alpha/\lambda}.$$

\smallskip
\noindent
Thus defined, $\bar{x}$ is an inflection point of $f$. Now if $f'(x_{2,0})\geq0$, then it should be ${\bar{x}} >x_{2,0}$ provided $\theta\in [0, b(\lambda,\alpha,p,x_{1,0})]$, but $a(\lambda,\alpha,p,x_{1,0},x_{2,0})< b(\lambda,\alpha,p,x_{1,0}),$ for the whole parametric set. Putting all together one can guarantee the existence of a unique root of $f$ in $(0,x_{2,0}]$ satisfying $f'(x)\geq 0$. See Figure \ref{fig:function} for an illustration of the behavior of $f$.
\end{obs}

\smallskip
\begin{theorem}
\label{T: WLLN}Consider the general Maki-Thompson model $\mathcal{X}_{GMT}(\theta,\lambda,\alpha,p;x_{1,0},x_{2,0},y_{1,0})$, and let $x_{1,\infty},x_{2,\infty}$ be given by Definition~
\ref{def:xinf2}.
Then,
\begin{equation*}
\lim_{N \to \infty} \, \frac{X_1^{(N)}(\tau^{(N)})}{N} = x_{1,\infty} \,\, \text{ and } \,\,
\lim_{N \to \infty} \, \frac{X_2^{(N)}(\tau^{(N)})}{N} = x_{2,\infty} \quad \text{in probability}. 
\end{equation*}
\end{theorem}

\smallskip
\begin{exa}
Consider the general Maki-Thompson model $\mathcal{X}_{GMT}(\theta,\lambda,\lambda,1;x_{1,0},x_{2,0},y_{1,0})$, i.e., let $\alpha=\lambda$ and $p=1$. We point out that this is the basic Maki-Thompson with directed inter-group interactions, where we assume groups of size $\theta\,N$ and $(1-\theta)\, N$, respectively. The resulting continuous-time Markov chain $\left\{\left(X^N_1(t),X^N_2(t),Y^N_1(t)\right)\right\}_{t\in[0,\infty)}$ has the following transitions and rates:
\begin{equation}\label{eq:rates_particular}
\begin{array}{cc}
\mbox{transition} \quad &\mbox{rate} \\[0.2cm]
(-1,0,1)       \quad & \lambda \, X_{1}\,Y_1, \\[0.2cm]
(0,-1,0)     \quad & \lambda\, Y_1\, X_2, \\[0.2cm]
(0,0,-1)      \quad & \lambda Y_1\,(N-X_1-X_2).\\[0.2cm]
\end{array}
\end{equation}
It is not difficult to see that $\lambda$ has no role in the behavior of the remaining proportion of ignorants, only in the speed of the process. However, with this choice of the parameters we can obtain $x_{1,\infty}$ and $x_{2,\infty}$ from Definition \ref{def:xinf2} explicitly in terms of the Lambert $W$ function, which is the multivalued inverse of the function $x \mapsto x\, e^x$. In other words, we obtain

$$x_{1,\infty}= -\frac{x_{1,0}}{2x_{1,0}+x_{2,0}}W_0\left(-(2x_{1,0}+x_{2,0})e^{-(y_{1,0}+2x_{1,0}+x_{2,0})}\right)$$ 
and
$$x_{2,{\infty}}=-\frac{x_{2,0}}{2x_{1,0}+x_{2,0}}W_0\left(-(2x_{1,0}+x_{2,0})e^{-(y_{1,0}+2x_{1,0}+x_{2,0})}\right),$$

\smallskip
\noindent
where $W_0$ is the principal branch of Lambert function $W$, see \cite{Corless} for more details about this function.
\end{exa}


The WLLN gives approximations for the proportions $X_1^{(N)}(\tau^{(N)})/N$ and $X_2^{(N)}(\tau^{(N)})/N$ for $N$ sufficiently large. Moreover, through stating a Central Limit Theorem, we can say how the fluctuations between these values and the asymptotical proportions $x_{1,\infty}$ and $x_{2,\infty}$, respectively, behave as $N$ growths.

\bigskip
\begin{theorem}\label{teoTCL}Consider the general Maki-Thompson model $\mathcal{X}_{GMT}(\theta,\lambda,\alpha,p;x_{1,0},x_{2,0},y_{1,0})$, and let $x_{1,\infty},x_{2,\infty}$ be given by Definition~
\ref{def:xinf2}. Then
	$$\sqrt{N}\left(\frac{X_1^{(N)}({\tau}^{(N)})}{N}-x_{1,{\infty}}, \frac{X_2^{(N)}({\tau}^{(N)})}{N}-x_{2,{\infty}}\right) \stackrel{\mathcal{L}}{\longrightarrow}  N_2(0,{\Sigma})~\mbox{ when }~N\rightarrow \infty,$$
where $ \stackrel{\mathcal{L}}{\longrightarrow} $ denotes convergence in law and $N(0,{\Sigma})$ is the bivariate normal distribution with mean zero and covariance matrix given by
\begin{equation*}
\Sigma =
\begin{pmatrix}
\Sigma_{11} & \Sigma_{12} \\
\Sigma_{21} & \Sigma_{22}
\end{pmatrix}
\end{equation*}

\smallskip
\noindent
whose elements may be written as explicit functions of the parameters and initial conditions, i.e. $\Sigma_{ij}:=\Sigma_{ij}\left(\theta, \lambda,\alpha,p,x_{1,0},x_{2,0},y_{1,0}\right)$ for $i,j\in \{1,2\}$
\end{theorem}

	\section{Proofs}\label{section:proofs}
	
Our proofs rely on a suitable application of well-known results from density dependent Markov chains. Roughly speaking, we consider a time-changed version of the model in such a way that the new process is a density dependent Markov chain. Thus defined, we can apply existing results about the convergence of these processes. This is a powerful technique which has been applied successfully before for other epidemic-like models, see \cite{arruda2,lebensztayn/machado/rodriguez/2011b,lebensztayn/machado/rodriguez/2011a}. In what follows consider the general Maki-Thompson model $\mathcal{X}_{GMT}(\theta,\lambda,\alpha,p;x_{1,0},x_{2,0},y_{1,0})$ and let ${V}^{(N)}:=(X_1^N(t),X_2^N(t),Y_1^N(t))$ for any $t\geq 0$.  	
	
	\subsection{The time-changed process and its deterministic limiting system}
	We consider a time-changed version of the original process, which is obtained by running the clock at rate $Y_1^{N}(t)^{-1}$. In other words, consider the continuous-time Markov chain $\{ \tilde{V}^{(N)}(t)\}_{t\in(0,\infty]}$, where $\tilde{V}^{(N)}(t)=(\tilde{X}_1^{N}(t),\tilde{X}_2^{N}(t),\tilde{Y}_1^{N}(t))$, with increments and corresponding rates given by
	
\begin{equation}
\begin{array}{ccc}
\mbox{transition} \quad &&\mbox{rate} \\[0.2cm]
l_1:=(-1,0,1)        && \lambda\, p\, {\tilde{X}}_{1},\\[0.2cm]
l_2:= (0,-1,0)      && \alpha\, {\tilde{X}}_2, \\[0.2cm]
l_3:=(0,0,-1)      && \alpha\, (N_2-{\tilde{X}}_2)+ \lambda\, (N_1-{\tilde{X}}_1),\\[0.2cm]
l_4:=(-1,0,0)       &&  \lambda\, (1-p)\,{\tilde{X}}_1.\\[.2cm]
\label{rates}
\end{array}
\end{equation}

\noindent
It is worth pointing out that this new Markov chain has the same transitions than the process $\{V^{(N)}(t)\}_{_{t\in(0,\infty]}}$ so if they start at the same point they should be absorbed also at the same point. In other words, if $\tilde{\tau}^{(N)}:=\inf\{t\geq 0: \tilde{Y}_1^{N}(t)=0\}$ what we have is
\begin{equation}
\label{eq:vv}
V^{(N)}\left(\tau^{(N)}\right) = \tilde{V}^{(N)}\left(\tilde{\tau}^{(N)}\right).
\end{equation}

\noindent
Now, let us consider the functions
$$
\begin{array}{rcl}
\beta_{l_1}(x_1,x_2,y_1) &:= &\lambda \,p\, x_1,\\[.2cm]
\beta_{l_2}(x_1,x_2,y_1) &:= &\alpha \,x_2,\\[.2cm]
 \beta_{l_3}(x_1,x_2,y_1)& := &\alpha\, (1-\theta -x_2) + \lambda\, (\theta-x_1),\\[.2cm]
 \beta_{l_4}(x_1,x_2,y_1)& :=&\lambda \,(1-p)\, x_1,
	\end{array}
	$$
\noindent
and note that the rates in \eqref{rates} can be written as $N{\beta}_{l_i}(\tilde{X}_1/N,\tilde{X}_2/N,\tilde{Y}_1/N)$, for $i\in \{1,2,3,4\}$. Thus defined, $\{{\tilde{V}}^{(N)}(t)\}_{t\in[0,\infty)}$ is a density dependent Markov chain with possible transitions in the set $\{l_1,l_2,l_3,l_4\}$. By \cite[Theorem 2.1 on Chapter 11]{EK} we conclude that the process $\{{\tilde{V}}^N (t)/N\}_{t\in[0,\infty)}$ converges almost surely as $N\to \infty$ to a deterministic limit determined by the drift function
\begin{equation}\label{eq:driftfunction}
F(x_1,x_2,y_1):=\displaystyle \sum_{i=1}^4 l_i {\beta}_{l_i}(x_1,x_2,y_1) =\left(-\lambda\, x_1, -\alpha \, x_2, \lambda\, [(1+p)\, x_1-\theta]-\alpha\, (1-\theta -x_2) \right).
\end{equation}

\noindent
More precisely, the deterministic system is governed by the following set of ordinary differential equations:
\begin{equation}\label{eq:systemEDO}
\left\{
\begin{array}{rcl}
x'_1(t)&=&-\lambda x_1,\\[.2cm]
x'_2(t)&=&-\alpha x_2,\\[.2cm]
y'_1(t)&=&\lambda[(1+p)x_1-\theta]-\alpha(1-\theta -x_2),\\[.2cm]
\end{array}
\right.
\end{equation}

\smallskip
\noindent
with initial conditions  $\left(x_1(0),x_2(0),y_1(0) \right)= (x_{1,0},x_{2,0},y_{1,0})$. The solution of \eqref{eq:systemEDO} is given by 

\begin{equation}\label{eq:system}
\left\{
\begin{array}{rcl}
x_1(t)&=&x_{1,0}e^{-\lambda t},\\[.2cm]
x_2(t)&=&x_{2,0}e^{-\alpha t},\\[.2cm]
y_1(t)&=&y_{1,0} +(1+p)(x_{1,0}-x_1(t))+\left(x_{2,0}-x_2(t)\right)+[\theta (\alpha-\lambda)-\alpha]t,\\[.2cm]
\end{array}
\right.
\end{equation}

\smallskip
\noindent
and we notice that after a straight calculation we have $y_1(t) = f(x_2(t))$, where $f$ is the function defined by Definition \ref{def:xinf2}. Now \cite[Theorem 2.1 on Chapter 11]{EK} guarantees that, on a suitable probability space,
\begin{equation}
\label{eq:ascv}
\displaystyle \lim_{N \to \infty} \frac{{\tilde{V}}^{(N)}(t)}{N}=v(t),~\mbox{a. s.}
\end{equation}
uniformly on bounded time intervals, where $v(t):=(x_1(t),x_2(t),y(t))$ for any $t\geq 0$. In particular, it can be proved, see \cite[Lemma 3.2]{lebensztayn/machado/rodriguez/2011b}, that
\begin{equation}
\label{eq:asc}
\displaystyle \lim_{N \to \infty} \frac{1}{N}\left({\tilde{X}}_1^{(N)}(t),{\tilde{X}}_2^{(N)}(t)\right)=(x_1(t),x_2(t)),~\mbox{a. s.}
\end{equation}
uniformly on $\mathbb{R}$.

\subsection{Proof of Theorems \ref{T: WLLN} and \ref{teoTCL}}

Both theorems follow from \cite[Theorem 4.1 on Chapter 11]{EK}. We adopt the notation used there, and refer the reader also to \cite{lebensztayn/machado/rodriguez/2011b} for an analogous proof of limit theorems to a related model. Although the technique of proof is standar, see \cite{lebensztayn/machado/rodriguez/2011a,lebensztayn/machado/rodriguez/2011b}, we summarize the main steps for the sake of completeness. Let $\varphi(x_1,x_2,y_1)=y_1$, and let ${\tau}_{\infty}:=\inf\{t:y_1(t)\leq 0\}$; then
\begin{equation}
\label{tau}
{\tau}_{\infty} = -\left(\frac{1}{\alpha}\right) \ln \left(\frac{x_{2,\infty}}{x_{2,0}}\right).
\end{equation}
Moreover, note that
\begin{equation}
\label{eq}
\varphi(v({\tau}_{\infty}))\cdot F(v({\tau}_{\infty}))={y_1}'({\tau}_{\infty})<0;
\end{equation}
which is obtained by the definition of $\tau_{\infty}$ and the fact that $y_1(t)$ is a concave function.

To prove Theorem \ref{T: WLLN} observe that that $y_{1,0}>0$ and (\ref{eq}) imply $y_1({\tau}_{\infty}-\epsilon)>0$ and $y_1({\tau}_{\infty}+\epsilon)<0$ for $0<\epsilon<{\tau}_{\infty}$. Then, from (\ref{eq:ascv}), follows that ${\tilde{Y_1}^{(N)}(t)}/{N}$ converges to $y_1$ almost surely uniformly on bounded  intervals, and hence 
\begin{equation}
\label{eq:tau}
\displaystyle \lim_{N \to \infty} {\tilde{\tau}}^{N}={\tau}_{\infty},~\mbox{a. s.}
\end{equation}
On the other hand, $y_{1,0}=0$, $y'_1(0)<0$ for all $t>0$, and the almost surely convergence of ${\tilde{Y_1}^{(N)}(t)}/{N}$ to ${y}_1$ uniformly on bounded intervals implies that \begin{equation}\label{eq:tau0}
\displaystyle \lim_{N \to \infty} {\tilde{\tau}}^{N}={\tau}_{\infty}=0,~\mbox{a. s.}
\end{equation}
By \eqref{eq:vv}, \eqref{eq:asc} and \eqref{eq:tau}/\eqref{eq:tau0} we obtain Theorem \ref{T: WLLN}. In order to prove Theorem \ref{teoTCL} we note from \cite[Teorema 4.1 on Chapter 11]{EK} that
$$\sqrt{N}\left(\frac{\tilde{X}_1\left(\tilde{{\tau}}^{(N)}\right)}{N} -x_{1,{\infty}}, \frac{\tilde{X}_2\left(\tilde{{\tau}}^{(N)}\right)}{N} -x_{2,{\infty}}\right)$$
converges in law, as $N \to \infty$, to 

\begin{equation}
\left(U_{x_1}({\tau}_{\infty}) +\left\{\frac{\lambda x_{1,{\infty}}}{{y_1}'({\tau}_{\infty})} \right\}U_{y_1} ({\tau}_{\infty})~,~ U_{x_2}({\tau}_{\infty})+\left\{\frac{\alpha x_{2,{\infty}}}{{y_1}'({\tau}_{\infty})}\right\}U_{y_1}({\tau}_{\infty})\right),
\label{tcl}
\end{equation}

\smallskip
\noindent
where $\mathcal{U}:=\left(U_{x_1},U_{x_2},U_{y_1}\right)$ is a multivariate Gaussian process with parameters to be defined later. The asymptotic distribution is a mean zero bivariate normal distribution. In what follows we explain how the covariance matrix $\Sigma$ is obtained. First, we find the matrix of partial derivatives of the drift function $F$ given by \eqref{eq:driftfunction} and the matrix $G$ defined in \cite[p. 458]{EK} by $G(x_1,x_2,y_1)=\sum_{i=1}^4 l_{i} l_{i}^T {\beta}_{l_i} (x_1,x_2,y_1).$ Therefore, we have

$$\partial F(x_1,x_2,y_1)=
\left(
\begin{array}{ccc}
-\lambda&0&0 \\
0&-\alpha&0\\
\lambda (1+p)&\alpha&0
\end{array}
\right)
$$

and

$$G(x_1,x_2,y_1)=
\left(
\begin{array}{ccc}
\lambda x_1&0&-\lambda p x_1 \\
0&\alpha x_2&0\\
-\lambda p x_1&0&\lambda(p-1)x_1-\alpha x_2 +(\lambda -\alpha)\theta -\alpha
\end{array}
\right).
$$

The next step is to obtain the solution $\phi$ of the matrix equation
$$\frac{\partial}{\partial t}\phi(t,s)=\partial F(x_1(t),x_2(t),y_1(t))\phi(t,s),~\phi(s,s)={\mathbb{I}}_3,$$

which is given by
$$\phi(t,s)=
\left(
\begin{array}{ccc}
e^{-\lambda(t-s)}&0&0 \\
0&e^{-\alpha(t-s)}&0\\
(1+p)\left(1-e^{-\lambda(t-s)}\right)&1-e^{-\alpha(t-s)}&1
\end{array}
\right).
$$

Finally, the covariance matrix of the Gaussian process $U$ at time $t$ is computed from
 
$$
Cov(U(t),U(t))\,=\,\int_0^t \phi(t,s)G(x_1(s),x_2(s),y_1(s))[\phi(t,s)]^T ds \label{cov}.
$$

\smallskip
\noindent
That is,
\begin{eqnarray}
Cov(U(t),U(t)) &=&\left(
\begin{array}{ccc}
c_{11}(t)& 0  & c_{13}(t)\\[.2cm]
0&c_{22}(t)&c_{23}(t)\\[.2cm]
c_{31}(t) & c_{32}(t) & c_{33}(t)
\end{array}
\right),
\end{eqnarray}

\noindent
where

\begin{eqnarray}
\label{c}
c_{11}(t)&=&\left(\frac{x_1(t)}{x_{1,0}}\right)(x_{1,0}-x_1(t)), \nonumber \\[.2cm]
c_{13}(t)&=&c_{31}(t)\;\;=\;\;x_1(t)\left[ \lambda t +(1+p)\left(1-\frac{x_1(t)}{x_{1,0}} \right) \right], \nonumber \\[.2cm]
c_{22}(t)&=& \left(\frac{x_2(t)}{x_{2,0}}\right)(x_{2,0}-x_2(t)), \\[.2cm]
c_{23}(t)&=&c_{32}(t)\;\;=\;\;x_2(t)\left(\alpha t-1+e^{-\alpha t}\right),\nonumber  \\[.2cm]
c_{33}(t)&=&x_2(t)\,(1+2\alpha -e^{-\alpha t} )-p\,(1+p)\,x_1(t)\,(e^{\lambda t}-\lambda t -1)+(p-1)\,(x_{1,0}-x_1(t))\nonumber \\[.2cm]
&& + \,t\, [(\lambda - \alpha)\theta+\alpha].\nonumber 
 \end{eqnarray}

Finally, we get the expression for $Cov(U({\tau}_{\infty}),U({\tau}_{\infty}))$ by observing that
$$e^{-\lambda \tau_{\infty}}=\frac{x_{1,{\infty}}}{x_{1,0}} \;\;\;\text{ and }\;\;\; \tau_{\infty} = -\left(\frac{1}{\alpha}\right) \ln \left(\frac{x_{2,\infty}}{x_{2,0}}\right),$$

\smallskip
\noindent
see \eqref{eq:system} and \eqref{tau}, respectively.

 
 \section{Aknowledgements}
 The second author thanks to Nicholas DiFonzo to pointing out the existence of the GBN-model to describe rumor spreading. This study was financed in part by the Coordena\c{c}\~ao de Aperfei\c{c}oamento de Pessoal de N\'ivel Superior - Brasil (CAPES) - Finance Code 001. This research was supported also by FAPESP (Grants 2016/11648-0, 2017/10555-0), and CNPq (Grant 304676/2016-0).

	\end{document}